\documentstyle[12pt,twoside]{article}

\font\teneufm=eufm10 scaled \magstep1
\font\seveneufm=eufm7 scaled \magstep1
\font\fiveeufm=eufm5  scaled \magstep1
\newfam\eufmfam
\textfont\eufmfam=\teneufm
\scriptfont\eufmfam=\seveneufm
\scriptscriptfont\eufmfam=\fiveeufm
\def\frak#1{{\fam\eufmfam\relax#1}}

\newfam\msbfam
\font\tenmsb=msbm10 scaled \magstep1  \textfont\msbfam=\tenmsb
\font\sevenmsb=msbm7 scaled \magstep1 \scriptfont\msbfam=\sevenmsb
\font\fivemsb=msbm5 scaled \magstep1  \scriptscriptfont\msbfam=\fivemsb
\def\Bbb{\fam\msbfam \tenmsb}

\def\RR{{\Bbb R}}
\def\CC{{\Bbb C}}

\def\PP{{\Bbb P}}

\def\ra{\rightarrow}

 \def\HollowBoxx #1#2#3{{\dimen0=#1 \advance\dimen0 by -#2       
       \dimen1=#1 \advance\dimen1 by #3                       
        \vrule height 0pt depth #3 width #2                   
       \hskip -#3
       \vrule height #1 depth #3 width #3}}                   
 \def\LeftContraction{\mathord{\kern1.45pt \HollowBoxx{6pt}{3.5pt}{.4pt}}\,}

 \def\HollowBox #1#2#3{{\dimen0=#1 \advance\dimen0 by -#3       
       \dimen1=#1 \advance\dimen1 by #3                       
        \vrule height #1 depth #3 width #3                    
        \vrule height 0pt depth #3 width #2                   
        \hskip -#3}}                                             
 \def\RightContraction{\mathord{\, \HollowBox{6pt}{3.1pt}{.4pt}} \kern1.6pt}             

\def\qed{{\hfill $\Box$}}
\newtheorem{theorem}{THEOREM}[section]
\newtheorem{lemma}[theorem]{Lemma}

\newtheorem{proposition}[theorem]{Proposition}

\begin{document}
\begin{center}
{\Large \bf On the Automorphism Groups 
\medskip\\
of Hyperbolic Manifolds}\footnote{{\bf Mathematics 
    Subject Classification:} 32H02, 32H20, 32M05}\footnote{{\bf Keywords 
   and Phrases:} hyperbolic complex manifolds, automorphism groups, equivalence problem}  
\medskip \\
\normalsize Alexander V. Isaev and Steven G. Krantz\footnote{Author
supported in part by a grant from the National Science Foundation.}
\end{center} 

\begin{quotation} \small \sl We show that there does not exist a
Kobayashi hyperbolic complex manifold of dimension $n\ne 3$,
whose group of holomorphic automorphisms
has dimension $n^2+1$ and that, if a 3-dimensional connected
hyperbolic complex manifold
has automorphism group of dimension 10, then it is
holomorphically equivalent to the Siegel space.  These results
complement earlier theorems of the authors on the possible
dimensions of automorphism groups of domains in complex space.

The paper also contains a proof of our earlier result
on characterizing $n$-dimensional hyperbolic complex
manifolds with automorphism groups of dimensions $\ge n^2+2$.  
\end{quotation}

\pagestyle{myheadings}
\markboth{A. V. Isaev and S. G. Krantz}{Automorphism Groups of
Hyperbolic Manifolds}

\setcounter{section}{0}

\section{Introduction}
\setcounter{equation}{0}

Let $M$ be a complex manifold of complex dimension $n$ and
$\hbox{Aut}(M)$ the group of its holomorphic automorphisms. The group
$\hbox{Aut}(M)$ is a topological group with the natural compact-open
topology. If $M$ is Kobayashi-hyperbolic (e.g., if $M$ is a bounded
domain in complex space), then $\hbox{Aut}(M)$ is in fact a
finite-dimensional real Lie group whose topology agrees with the compact-open topology
\cite{Ko1}, \cite{Ko2}.

We are interested in the problem of characterizing hyperbolic complex
manifolds by the dimensions of their automorphism groups. Let $M$ be
such a manifold. It is known (see \cite{Ka}, \cite{Ko1}) that
$\hbox{dim}\,\hbox{Aut}\,(M)\le n^2+2n$ and, if $M$ is connected and 
$\hbox{dim}\,\hbox{Aut}\,(M)=n^2+2n$, then $M$ is holomorphically
equivalent to the unit ball $B^n\subset\CC^n$. In \cite{IK} we
obtained the following result.

\begin{theorem}\label{char} \sl Let $M$ be a connected hyperbolic manifold of
  complex dimension $n\ge 2$. Then the following holds:
\smallskip\\

\noindent {\bf (a)} If $\hbox{dim}\,\hbox{Aut}(M)\ge n^2+3$, then $M$
is biholomorphically equivalent to $B^n$.
\smallskip\\

\noindent {\bf (b)}  If $\hbox{dim}\,\hbox{Aut}(M)= n^2+2$, then $M$
is biholomorphically equivalent to $B^{n-1}\times\Delta$, where
$\Delta$ is the unit disc in $\CC$.
\end{theorem}

In this paper we completely characterize hyperbolic manifolds
with automorphism groups of
dimension $n^2+1$.
In \cite{GIK} we observed that the automorphism group of a hyperbolic
Reinhardt domain in $\CC^n$ cannot have dimension $n^2+1$.
Therefore, it has been our expectation that there should be very few manifolds
with $\hbox{dim}\,\hbox{Aut}(M)=n^2+1$. In fact, we have known of
only one example of a manifold with such an automorphism group
dimension; this occurs in dimension $n=3$.
\medskip\\

\noindent {\bf Example.} Consider the 3-dimensional Siegel space
(the symmetric bounded domain of type $(III_2)$):
$$
S:=\left\{(z_1,z_2,z_3)\in\CC^3: E-Z\overline{Z}\,\,\hbox{is a
    positive-definite matrix}\right\},
$$
where
$$
Z:=\left(
\begin{array}{cc}
\displaystyle z_1 & \displaystyle z_2\\
\displaystyle z_2 & \displaystyle z_3
\end{array}
\right),
$$
and $E$ is the $2\times 2$ identity matrix. The automorphism group of this domain is isomorphic to
$Sp_4(\RR)/\{\pm \hbox{Id}\}$ and has dimension $10=3^2+1$ (see, e.g., \cite{S}).
\smallskip\\

In this paper we prove the following theorem that shows that the
Siegel space $S$ is indeed the only possibility.

\begin{theorem}\label{char2} \sl Let $M$ be a connected hyperbolic
  manifold of complex dimension $n\ge 2$. 
 Then the following holds:
\smallskip\\ 

\noindent {\bf (a)} If $n\ne 3$, then $\hbox{dim}\,\hbox{Aut}(M)\ne n^2+1$.
\smallskip\\

\noindent {\bf (b)} If $n=3$ and
$\hbox{dim}\,\hbox{Aut}(M)=10$, then $M$ is holomorphically equivalent
to the Siegel space $S$.
\end{theorem}

Note that, for $n=1$, all manifolds with positive-dimensional
automorphism groups are known explicitly \cite{FK}. In particular, if
$\hbox{dim}\,\hbox{Aut}(M)=3$, then $M$ is equivalent to $\Delta$, and
there does not exist a hyperbolic manifold $M$ with $\hbox{dim}\,\hbox{Aut}(M)=2$. 

There is probably no hope to obtain a complete classification of all $n$-dimensional
hyperbolic manifolds with automorphism groups of dimension $n^2$ or
smaller. Indeed, in \cite{GIK} we classified all hyperbolic Reinhardt
domains in $\CC^n$ whose automorphism groups have dimension $n^2$, and
even in this special case the classification is rather large and non-trivial.

For the sake of completeness of our exposition, and because
some of the ideas originating there will recur in
the proof of Theorem \ref{char2}, we reproduce the
proof of Theorem \ref{char} in Section 2. We prove Theorem \ref{char2}
in Section 3.

Before proceeding, we would like to thank W. Kaup for helpful
discussions and for showing us useful references.

\section{Proof of Theorem \ref{char}}
\setcounter{equation}{0}

For the proof of Part {\bf (a)} of Theorem 1.1 we need the following
lemma.

\begin{lemma} \label{alg} \sl Let $G$ be a Lie subgroup of the unitary
group $U(n)$ and let $G^c$ be its connected component of the identity.
Suppose that $\hbox{dim}\,G\ge n^2-2n+3$, $n\ge 2$, $n\ne 4$. Then either $G=U(n)$,
or $G^c=SU(n)$. For $n=4$ this list must be augmented by subgroups
of $U(4)$ whose Lie algebras are isomorphic to $\RR \oplus {\frak {sp}}_{2,0}$.  
\end{lemma}

\noindent{\bf Proof of Lemma \ref{alg}:} Since $G$ is compact, it is completely reducible
(see, e.g., \cite{VO}) and thus is isomorphic to a direct product
$G_1\times\dots\times G_k$; here $G_j$ for each $j$ is a compact
subgroup of $U(n_j)$, $\sum_{j=1}^k n_j=n$, and $G_j$ acts irreducibly
on $\CC^{n_j}$. Since $\hbox{dim}\,G_j\le n_j^2$ and $\hbox{dim}\,G\ge
n^2-2n+3$, it follows that $k=1$, i.e. $G$ acts (complex) irreducibly
on $\CC^n$.

Let ${\frak g}\subset{\frak {u}}_n$ be the Lie algebra of $G$ and
${\frak g}^{\CC}\subset{\frak {gl}}_n$ its complexification. It then
follows that ${\frak g}^{\CC}$ acts irreducibly and faithfully on $\CC^n$. Therefore
by a theorem of \'E. Cartan (see, e.g., \cite{GG}), ${\frak g}^{\CC}$ is
either semisimple or is the direct sum of a semisimple ideal ${\frak
h}$ and $\CC$, where $\CC$ acts on $\CC^n$ by multiplication. Clearly
the action of the ideal ${\frak h}$ on $\CC^n$ is irreducible and faithful.

Suppose first that ${\frak g}^{\CC}$ is semisimple, and let ${\frak
g}^{\CC}={\frak g}_1\oplus\dots\oplus{\frak g}_m$ be its decomposition
into the direct sum of simple ideals. It then follows (see,
e.g., \cite{GG}) that the representation of ${\frak g}^{\CC}$ is the
tensor product of some irreducible faithful representations of ${\frak
g}_j$. Let $n_j$ denote the dimension of the representation of ${\frak
g}_j$, $j=1,\dots,m$. Then $n=n_1\cdot\dots\cdot n_m$ and
$\hbox{dim}_{\CC}\,{\frak g}_j\le n_j^2-1$, $n_j\ge 2$ for
$j=1,\dots,m$. It is now not difficult to prove the following claim.

\begin{quote}
{\bf Claim:} {\sl If $n=n_1\cdot\dots\cdot n_m$, $m\ge 2$, $n_j\ge 2$
for $j=1,\dots,m$, then $\sum_{j=1}^m n_j^2\le n^2-2n$.}
\end{quote}

It follows from the claim that $m=1$, i.e., ${\frak g}^{\CC}$ is
simple. The minimal dimensions of irreducible faithful
representations of complex simple Lie algebras are well-known (see,
e.g., \cite{VO}). In the table below $V$ denotes representations of
minimal dimension.

\begin{center}
\begin{tabular}{|l|c|c|}
\hline
\multicolumn{1}{|c|}{${\frak g}$}&
\multicolumn{1}{c|}{ $\hbox{dim}\,V$}&
\multicolumn{1}{c|}{$\hbox{dim}\,{\frak g}$}
\\ \hline
${\frak {sl}}_k$\,\,$k\ge 2$ & $k$ & $k^2-1$
\\ \hline
${\frak o}_k$\,\, $k\ge 7$&  $k$ & $\frac{k(k-1)}{2}$   
\\ \hline
${\frak {sp}}_{2k}$\,\,$k\ge 2$ & $2k$ & $2k^2+k$ 
\\ \hline
${\frak e}_6$ & 27 & 78
\\ \hline
${\frak e}_7$ & 56 & 133
\\ \hline
${\frak e}_8$ & 248 & 248
\\ \hline
${\frak f}_4$ & 26 & 52
\\ \hline
${\frak g}_2$ & 7 & 14
\\ \hline
\end{tabular}
\end{center}

Since $\hbox{dim}_{\CC}\,{\frak g}^{\CC}\ge n^2-2n+3$, it follows that
${\frak g}^{\CC}\simeq{\frak {sl}}_n$. Since ${\frak g}$ is a compact
algebra, we get that ${\frak g}={\frak {su}}_n$ (see \cite{VO}) and
therefore $G^c=SU(n)$ (note here that if ${\frak g}$ is a subalgebra in
${\frak u}_n$ and ${\frak g}$ is isomorphic to ${\frak {su}}_n$, then
${\frak g}$ coincides with ${\frak {su}}_n$, i.e., it consists exactly of
matrices with zero trace).

Suppose now (for the second case) that ${\frak g}^{\CC}={\frak
h}\oplus\CC$, where ${\frak h}$ is a semisimple ideal in ${\frak
g}^{\CC}$. Then, repeating the above argument for ${\frak h}$ and
taking into account that $\hbox{dim}\,{\frak h}\ge n^2-2n+2$, we
conclude that ${\frak h}\simeq{\frak {sl}}_n$ for $n\ne 4$. 
Therefore, for $n\ne 4$, ${\frak g}^{\CC}={\frak {gl}}_n$ and hence
${\frak g}={\frak u}_n$, which implies that $G=U(n)$. 

For $n=4$, either ${\frak h}\simeq{\frak {sl}}_4$ or ${\frak
h}\simeq{\frak {sp}}_4$.  We thus find that either ${\frak g}={\frak
u}_4$ (in which case $G=U(4)$), or ${\frak g}\simeq\RR\oplus{\frak
{sp}}_{2,0}$.

The lemma is proved. \qed
\medskip  \\

\noindent {\bf Proof of Theorem \ref{char}, Part (a):} Let $p\in M$ and let 
$I_p$ denote the isotropy group of $p$ in $\hbox{Aut}(M)$. Since the
complex dimension of $M$ is $n$, the real dimension of any orbit of
the action of $\hbox{Aut}(M)$ on $M$ does not exceed $2n$, and
therefore we have $\hbox{dim}\,I_p\ge n^2-2n+3$.

Consider the isotropy representation $\alpha_p: I_p\ra GL(T_p(M),\CC)$:  
$$
\alpha_p(f):=df(p), \qquad f\in I_p.
$$
The mapping $\alpha_p$ is a continuous group homomorphism (see,
e.g., Lemma 1.1 of \cite{GK}) and thus is a Lie group homomorphism (see
\cite{Wa}). Since $I_p$ is compact (see \cite{Ko1}), there is a positive-definite Hermitian form
$h_p$ on $T_p(M)$ such that
$\alpha_p(I_p)\subset U_{h_p}(n)$, where $U_{h_p}(n)$ is the group of
complex linear
transformations of $T_p(M)$ preserving the form $h_p$. We choose a
basis in $T_p(M)$ such that $h_p$ in this basis is given by the
identity matrix.
 
By \cite{E} and \cite{Ki}, the mapping $\alpha_p$ is one-to-one.
Further, since $\hbox{dim}\,I_p\ge n^2-2n+3$, we see that $\alpha(I_p)$ is a compact subgroup of
$U_{h_p}(n)$ of dimension at least $n^2-2n+3$.  We are now going to
use Lemma \ref{alg}. 

Assume first that $n\ne 4$. Then  we have that either
$\alpha_p(I_p)=U_{h_p}(n)$, or $\alpha_p(I_p)^c=SU_{h_p}(n)$ (the
latter denotes the subgroup of $U_{h_p}(n)$ consisting of matrices
with determinant 1). The groups $U_{h_p}(n)$ and $SU_{h_p}(n)$ act
transitively on the unit sphere in $T_p(M)$ and thus act transitively
on directions in $T_p(M)$ (see \cite{GK} and \cite{BDK} for
terminology). Since $M$ is non-compact (because the dimension of
$\hbox{Aut}(M)$ is positive---see \cite{Ko1}), the main result of
\cite{GK} and its generalization in \cite{BDK} applies.  Thus $M$ is
biholomorphically equivalent to $B^n$ (and therefore the possibility
$\alpha_p(I_p)^c=SU_{h_p}(n)$ is in fact not realizable).

Suppose now that $n=4$. If we have that either
$\alpha_p(I_p)=U_{h_p}(4)$, or else
$\alpha_p(I_p)^c=SU_{h_p}(4)$ for some $p\in M$, then by the above argument $M$ is equivalent to
$B^4$. Suppose now that the Lie algebra of $\alpha_p(I_p)$ is
isomorphic to $\RR\oplus{\frak {sp}}_{2,0}$ for every $p\in M$. Then
$\hbox{dim}\,\alpha_p(I_p)=11$ for any $p$. Since
$\hbox{dim}\,\hbox{Aut}(M)\ge 19$, we have that in fact 
$\hbox{dim}\,\hbox{Aut}(M)=19$, and thus $M$ is homogeneous. Therefore, by
\cite{N}, \cite{P-S}, $M$ is biholomorphically equivalent to a Siegel domain
$D\subset\CC^4$ of the first or second kind. Further, we note that any
representation $\phi:{\frak {sp}}_{2,0}\ra{\frak {gl}}_4$ is conjugate
to the standard embedding of ${\frak {sp}}_{2,0}$ into ${\frak {gl}}_4$ by an element from $GL(4,\CC)$ (to see this, one can
extend
$\phi$ to a 4-dimensional representation of the complex
Lie algebra ${\frak {sp}}_4$ and notice that such a representation
is unique up to conjugation by elements of $GL(4,\CC)$).   
Therefore $\phi({\frak
  {sp}}_{2,0})$ contains an element $X$ such that
$\exp(X)=-\hbox{id}$, and thus $\alpha_p(I_p)$ contains $-\hbox{id}$ for any
$p$.  Hence the domain $D$ is in fact symmetric. It now follows from
the explicit classification of symmetric Siegel domains (see
\cite{S}) that in fact there is no symmetric Siegel domain in $\CC^4$ with
automorphism group of dimension equal to 19. This concludes the proof of Part {\bf (a)}.

For the proof of Part {\bf (b)} we will need the following lemma (which
follows from the proof of Lemma \ref{alg}).

\begin{lemma}\label{alg1} \sl Let $U_h(n)$ be the group of linear
  transformations of a complex $n$-dimensional space $V$ that preserve
  a positive-definite Hermitian form $h$ on $V$, and let $G$ be a Lie
  subgroup of $U_h(n)$ with $\hbox{dim}\,G\ge n^2-2n+2$, $n\ge 2$, $n\ne
4$. Then either $G=U_h(n)$, or $G^c=SU_h(n)$, or $V$ splits into a sum
  of 1- and $(n-1)$-dimensional $h$-orthogonal complex subspaces $V^1$ and
  $V^2$ such that $G=U_{h^1}(1)\times U_{h^2}(n-1)$, where $h^j$ is
  the restriction of $h$ to $V^j$. For $n=4$, there is the additional
  possibility that $G$ can be any
  subgroup of $U_h(4)$ with Lie algebra isomorphic to either 
${\frak {sp}}_{2,0}$ or $\RR\oplus{\frak {sp}}_{2,0}$.
\end{lemma}

\noindent {\bf Proof of Theorem \ref{char}, Part (b):} We will use the notation
from the proof of Part {\bf (a)} above. Let $p\in M$ and $I_p$ be the isotropy group
of $p$ in $\hbox{Aut}(M)$. Then we have $\hbox{dim}\,I_p\ge
n^2-2n+2$ and thus $\alpha_p(I_p)$ is a subgroup of $U_{h_p}$ of dimension at
least $n^2-2n+2$. We now use Lemma \ref{alg1}. If, for some $p\in M$,
we have that either $\alpha_p(I_p)=U_{h_p}(n)$ or
$\alpha_p(I_p)^c=SU_{h_p}(n)$, then $\alpha_p(I_p)$ acts
transitively on directions in $T_p(M)$.   Hence, as
in the proof of Part {\bf (a)}, $M$ is biholomorphically
equivalent to $B^n$; this is impossible since
$\hbox{dim}\,\hbox{Aut}(M)=n^2+2$. 

Further suppose that, for any point $p\in M$, $T_p(M)$ splits into
the sum of 1- and $(n-1)$-dimensional $h_p$-orthogonal complex subspaces $V_p^1$ and
$V_p^2$ such that $\alpha_p(I_p)=U_{h_p^1}(1)\times U_{h_p^2}(n-1)$. In particular,
$\hbox{dim}\,I_p=n^2-2n+2$ for all $p\in M$ and
therefore $M$ is homogeneous. Then, by \cite{N}, \cite{P-S}, $M$ is
biholomorphically equivalent to a homogeneous Siegel domain $D$ of the
first or second kind. Since $\alpha_p(I_p)$ contains the
transformation $-\hbox{id}$ for
all $p\in M$, the domain $D$ is in fact symmetric. The theorem for
$n\ne 4$ now
follows from the explicit classification of symmetric Siegel domains
(see \cite{S}).

Suppose now that $n=4$ and that, for some point $p\in M$, the Lie
algebra of $\alpha_p(I_p)$ is isomorphic to either ${\frak
{sp}}_{2,0}$ or to $\RR\oplus{\frak {sp}}_{2,0}$. In the proof of
Part {\bf (a)} we noted that any embedding of ${\frak
{sp}}_{2,0}$ into ${\frak {gl}}_4$ is conjugate by an element of
$GL(4,\CC)$ to the standard one. Therefore $\alpha_p(I_p)$ contains a
subgroup that is conjugate by an element of $GL(4,\CC)$ to $Sp_{2,0}$. Since
$Sp_{2,0}$ acts transitively on the sphere of dimension 7, we get
that $\alpha_p(I_p)$ acts transitively on directions in $T_p(M)$ and
therefore, as in the proof of Part {\bf (a)}, $M$ is
biholomorphically equivalent to the unit ball which is impossible.

The theorem is proved. \qed
\medskip\\     

\noindent{\bf Remark.} The argument in the last paragraph in the proof
of Part {\bf (b)} could also be used in the proof of Part {\bf (a)} for the case $n=4$ 
without reference to the classification
theory of symmetric domains.  For hyperbolic Reinhardt domains,
Theorem \ref{char} was obtained by a different argument in
\cite{GIK}.

\section{Proof of Theorem \ref{char2}}
\setcounter{equation}{0}

We will use the notation from the proof of Theorem \ref{char}. Suppose
first that $n=2$. Since $\hbox{dim}\,\hbox{Aut}(M)=5$, for any point $p\in M$ we have $\hbox{dim}\,\alpha_p(I_p)\ge 1$.
It follows from Lemma \ref{alg1} that, if $G$ is a positive-dimensional
closed subgroup of $U(2)$, then one of the following holds:
\medskip\\

\noindent {\bf (i)} $G=U(2)$;
\smallskip\\

\noindent {\bf (ii)} $G^c=SU(2)$;
\smallskip\\

\noindent {\bf (iii)} $\CC^2$ splits into a sum of two 1-dimensional
orthogonal subspaces $V^1$ and $V^2$, such that $G=U_{h^1}(1)\times
U_{h^2}(1)$, where $h^j$ is the restriction of the standard Hermitian
form on $\CC^2$ to $V^j$;
\smallskip\\

\noindent {\bf (iv)} $G$ is 1-dimensional.
\medskip\\

If, for some $p\in M$, $\alpha_p(I_p)$ is as in {\bf (i)} or {\bf
  (ii)}, then $M$ is holomorphically equivalent to
  $B^2\subset\CC^2$, which is impossible since
  $\hbox{dim}\,\hbox{Aut}(M)=5$. 

If, for every $p\in M$,
  $\alpha_p(I_p)$ is as the group $G$ in {\bf (iii)} (i.e., $T_p(M)$
  splits into a sum of subspaces $V^1$ and $V^2$ as in {\bf (iii)}), 
then, for every $p\in M$,
  $\alpha_p(I_p)$ contains the element $-\hbox{id}$. Since $M$ is
  hyperbolic, it is a complex metric Banach manifold (when equipped
  with the Kobayashi metric). It now follows from the results of \cite{V}
  (see Theorem 17.16 in \cite{U}) that $M$ is homogeneous, which is
  again impossible. 

If, for every $p\in M$, $\alpha_p(I_p)$ is 1-dimensional, then $M$ is
  homogeneous and hence is equivalent to a bounded homogeneous domain
  in $\CC^2$ \cite{N}. Therefore, $M$ is equivalent to either $B^2$ or
  $\Delta^2$. But such an equivalence is impossible since
  $\hbox{dim}\,\hbox{Aut}(B^2)=8$ and $\hbox{dim}\,\hbox{Aut}(\Delta^2)=6$.

We now assume that $M=M_1\cup M_2$, with $M_1, M_2\ne \emptyset$, and, for
$p\in M_1$, $\alpha_p(I_p)$ is as the group $G$ in {\bf (iii)}
and, for $p\in M_2$,
$\alpha_p(I_p)$ is 1-dimensional. Fix $p_0\in M_2$ and let $\Omega$ be
the orbit of $p_0$ under $\hbox{Aut}(M)^c$. Then $\Omega$ is a
homogeneous subdomain of $M$. Therefore, $\Omega$ is
holomorphically equivalent to either $B^2$ or $\Delta^2$.
Consider the restriction map $\phi: \hbox{Aut}(M)^c\ra \hbox{Aut}(\Omega)$:
$$
\phi: f\mapsto f|_{\Omega}.
$$
Clearly, $\phi$ is continuous and hence a Lie group homomorphism. It is
also one-to-one by the uniqueness theorem.  
Therefore
$\hbox{Aut}(\Omega)$ contains a (not necessarily closed) subgroup $H$ of dimension
5 that acts transitively on $\Omega$. Since $M_1\ne\emptyset$, $H$
contains a subgroup isomorphic to $U(1)\times U(1)$. It now follows from the explicit formulas for the
automorphism groups of $B^2$ and $\Delta^2$ that such a
subgroup $H$ in fact does not exist. 
This proves Part {\bf (a)} for $n=2$.

Suppose now that $n=3$. If $\hbox{dim}\,\hbox{Aut}(M)=10$ then, for any
$p\in M$, we have $\hbox{dim}\,\alpha_p(I_p)\ge 4$.
It follows from Lemma \ref{alg1} that, if $G$ is a
closed subgroup of $U(3)$ of dimension at least 4, then one of the
following holds:
\medskip\\

\noindent {\bf (i)} $G=U(3)$;
\smallskip\\

\noindent {\bf (ii)} $G^c=SU(3)$;
\smallskip\\

\noindent {\bf (iii)} $\CC^3$ splits into a sum of a 1- and
2-dimensional orthogonal subspaces $V^1$ and $V^2$ respectively such that
$G=U_{h^1}(1)\times U_{h^2}(2)$, where $h^j$ is the restriction of the
standard Hermitian form on $\CC^3$ to $V^j$;
\smallskip\\

\noindent {\bf (iv)} $G$ is 4-dimensional.
\medskip\\

If, for some $p\in M$, $\alpha_p(I_p)$ is as in {\bf (i)} or {\bf
  (ii)} then, as above, $M$ is holomorphically equivalent to
  $B^3\subset\CC^3$, which is impossible since
  $\hbox{dim}\,\hbox{Aut}(M)=10$. 

If, for every $p\in M$,
  $\alpha_p(I_p)$ is as the group $G$ in {\bf (iii)}, then, by
  \cite{V}, as in the case $n=2$ above, $M$ is homogeneous which is impossible. 

If, for every
  $p\in M$, $\alpha_p(I_p)$ is 4-dimensional, then $M$ is homogeneous
  and hence is holomorphically equivalent to one of the following
  domains: $B^3$, $B^2\times\Delta$, $\Delta^3$, or the Siegel space
$S$. Among these domains only $S$ has automorphism group of dimension 10. 

We now assume that $M=M_1\cup M_2$, with $M_1, M_2\ne \emptyset$, and, for $p\in M_1$, we suppose that $\alpha_p(I_p)$ is as the group $G$ in {\bf
(iii)} and, for $p\in M_2$, $\alpha_p(I_p)$ is 4-dimensional.
As in the case $n=2$ above, this implies that there exists a
subdomain
$\Omega\subset M$ such that $\hbox{Aut}(\Omega)$ contains a
subgroup $H$ of dimension 10 that acts transitively on $\Omega$ and
that contains a subgroup isomorphic to $U(1)\times U(2)$. It now
follows from the explicit formulas for the automorphism groups of
$B^3$, $B^2\times\Delta$ and $S$ that such a subgroup $H$ in fact does
not exist. This proves Part {\bf (b)}.

Let now $n\ge 4$. Since $\hbox{dim}\,\hbox{Aut}(M)=n^2+1$,
for any $p\in M$ we have $\hbox{dim}\,\alpha_p(I_p)\ge n^2-2n+1$. It
follows from Lemma \ref{alg1} that, if $G$ is a closed
subgroup of $U(n)$ of dimension at least $n^2-2n+1$, then one of the
following holds:
\medskip\\

\noindent {\bf (i)} $G=U(n)$;
\smallskip\\

\noindent {\bf (ii)} $G^c=SU(n)$;
\smallskip\\

\noindent {\bf (iii)} $\CC^n$ splits into a sum of a 1- and
$(n-1)$-dimensional orthogonal subspaces $V^1$ and $V^2$ respectively such that
$G=U_{h^1}(1)\times U_{h^2}(n-1)$, where $h^j$ is the restriction of the
standard Hermitian form on $\CC^n$ to $V^j$;
\smallskip\\

\noindent {\bf (iv)} $\hbox{dim}\,G=n^2-2n+1$;
\medskip\\

and, for $n=4$, there is one more possibility:
\medskip\\

\noindent {\bf (v)} The Lie algebra of $G$ is isomorphic to either
${\frak {sp}}_{2,0}$ or $\RR\oplus{\frak {sp}}_{2,0}$.
\medskip\\

If, for some $p\in M$, $\alpha_p(I_p)$ is as in {\bf (i)}, {\bf
  (ii)} or {\bf (v)}, then, as above, $M$ is holomorphically equivalent to
  $B^n\subset\CC^n$, which is impossible since
  $\hbox{dim}\,\hbox{Aut}(M)=n^2+1$. 

If, for every $p\in M$,
  $\alpha_p(I_p)$ is as the group $G$ in {\bf (iii)}, then by
  \cite{V}, $M$ is homogeneous which is impossible. 

If, for every
  $p\in M$, $\alpha_p(I_p)$ is $n^2-2n+1$-dimensional, then $M$ is
  homogeneous and hence, by \cite{N}, \cite{P-S}, is equivalent to a Siegel domain of the first
  or second kind in $\CC^n$. We now need the
  following proposition which is of independent interest as it gives
  an alternative proof of Theorem \ref{char} and a proof of Theorem
  \ref{char2} for Siegel domains in $\CC^n$, $n\ge 4$.

\begin{proposition}\label{Sig} \sl Let $U\subset\CC^n$, $n\ge 4$, be a
  Siegel domain of the first or second kind. Suppose that
  $\hbox{dim}\,\hbox{Aut}(U)\ge n^2+1$. Then $U$ is
  holomorphically equivalent to either $B^n$ or
  $B^{n-1}\times\Delta$.
\end{proposition}

\noindent{\bf Proof of Proposition \ref{Sig}:} The domain $U$ has the
  form
\begin{equation}
U=\left\{(z,w)\in\CC^{n-k}\times\CC^k:\hbox{Im}\,w-F(z,z)\in
    C\right\}, \label{om}
\end{equation}
where $1\le k\le n$, $C$ is an open convex cone in $\RR^k$ not
containing an entire line and $F=(F_1,\dots,F_k)$ is a $\CC^k$-valued Hermitian form on
$\CC^{n-k}\times\CC^{n-k}$ such that
$F(z,z)\in\overline{C}\setminus\{0\}$ for all non-zero $z\in\CC^{n-k}$.

We will first show that $k\le 2$. It follows from \cite{KMO} that
$$
\hbox{dim}\,\hbox{Aut}(U)\le 4n-2k+\hbox{dim}\,{\frak g}_0(U).
$$
Here ${\frak g}_0(U)$ is the Lie algebra of all vector fields on $\CC^n$
of the form
$$
X_{A,B}=Az\frac{\partial}{\partial z}+Bw\frac{\partial}{\partial w},
$$
where $A\in{\frak {gl}}_{n-k}$, $B$ belongs to the Lie algebra ${\frak
  g}(C)$ of the
group $G(C)$ of linear automorphisms of the cone $C$, and the following holds:
\begin{equation}
F(Az,z)+F(z,Az)=BF(z,z), \label{iden}
\end{equation}
for $z\in\CC^{n-k}$.
By the definition of Siegel domain, there exists a positive-definite
linear combination $R$ of the components of the Hermitian form
$F$. Then, for a fixed matrix $B$ in formula (\ref{iden}), the matrix
$A$ is determined at most up to a matrix that is Hermitian with
respect to $R$. Since the dimension of the algebra of matrices
Hermitian with respect to $R$ is equal to $(n-k)^2$, we have
$$
\hbox{dim}\,{\frak g}_0(U)\le (n-k)^2+\hbox{dim}\,{\frak
  g}(C),
$$
and thus the following holds
\begin{equation}
\hbox{dim}\,\hbox{Aut}(U)\le 4n-2k+(n-k)^2+\hbox{dim}\,{\frak
  g}(C). \label{d}
\end{equation}

\begin{lemma}\label{est} \sl We have 
$$
\hbox{dim}\,{\frak g}(C)\le \frac{k^2}{2}-\frac{k}{2}+1.
$$
\end{lemma}

\noindent{\bf Proof of Lemma \ref{est}:} Fix a point $x_0\in C$
and consider its stabilizer $G_{x_0}(C)\subset G(C)$.  The stabilizer is
compact since it leaves stable the bounded open set $C\cap (x_0-C)$ and
therefore we can assume that it is contained in the group
$O(k,\RR)$. The group $O(k,\RR)$ acts transitively on the sphere $S(|x_0|)$ of
radius $|x_0|$ in $\RR^k$, and the stabilizer $H_{x_0}$ of the point
$x_0\in S(|x_0|)$ under this action is isomorphic to
$O(k-1,\RR)$. Since $G_{x_0}\subset H_{x_0}$, we have
$$
\hbox{dim}\,G_{x_0}\le
\hbox{dim}\,H_{x_0}=\frac{k^2}{2}-\frac{3k}{2}+1,
$$
and therefore
$$
\hbox{dim}\,{\frak g}(C)\le \frac{k^2}{2}-\frac{k}{2}+1.
$$

The lemma is proved. \qed
\medskip\\

It now follows from (\ref{d}) and Lemma \ref{est} that
\begin{equation}
\hbox{dim}\,\hbox{Aut}(U)\le
\frac{3k^2}{2}-k\left(2n+\frac{5}{2}\right)+n^2+4n+1. \label{d1}
\end{equation}
It is easy to check that the right-hand side in (\ref{d1}) is strictly
less than $n^2+1$ for $n\ge 4$ and $k\ge 3$. Therefore, $k\le 2$.

If $k=1$, the domain $U$ is equivalent to $B^n$. Suppose that
$k=2$. Without loss of generality we can assume that the first
component $F_1$ of the $\CC^2$-valued Hermitian form $F$, is
positive-definite. We will show that the second component $F_2$ has to be
proportional to $F_1$. Indeed, if $F_2$ is not proportional to $F_1$,
then in formula (\ref{iden}) the matrix $A$ is determined by the
matrix $B$ up to transformations that are : (1) Hermitian with respect
to the positive-definite form $F_1$; (2) Hermitian with respect to
some other Hermitian form which is not proportional to $F_1$. The dimension of the algebra of
matrices satisfying conditions (1) and (2) does not exceed
$(n-2)^2-2$, and therefore
$$
\hbox{dim}\,{\frak g}_0(U)\le (n-2)^2-2+\hbox{dim}\,{\frak
  g}(C).
$$
Hence
$$
\hbox{dim}\,\hbox{Aut}(U)\le n^2-2+\hbox{dim}\,{\frak
  g}(C),
$$
which together with Lemma \ref{est} implies
$$
\hbox{dim}\,\hbox{Aut}(U)\le n^2,
$$
which is a contradiction. Thus, $F_2$ is proportional to
$F_1$. Therefore, $U$ is holomorphically equivalent to one of the
following domains:
$$
U_1:=\left\{(z,w)\in\CC^{n-2}\times\CC^2: \hbox{Im}\,w_1-|z|^2>0,\,\hbox{Im}\,w_2>0\right\},
$$
or
$$
U_2:=\left\{(z,w)\in\CC^{n-2}\times\CC^2: \hbox{Im}\,w_1-|z|^2>0,\,\hbox{Im}\,w_2-|z|^2>0\right\}.
$$

The domain $U_1$ is equivalent to $B^{n-1}\times\Delta$. We will show
that $\hbox{dim}\,\hbox{Aut}(U_2)<n^2+1$. Let ${\frak g}(U_2)$ be the Lie algebra of
$\hbox{Aut}(U_2)$. By \cite{KMO}, ${\frak g}(U_2)$ is a graded Lie
algebra:
$$
{\frak g}(U_2)={\frak g}_{-1}(U_2)\oplus {\frak g}_{-1/2}(U_2)\oplus
{\frak g}_0(U_2)\oplus {\frak g}_{1/2}(U_2)\oplus {\frak g}_1(U_2),
$$
where $\hbox{dim}\,{\frak g}_{-1}(U_2)=2$, $\hbox{dim}\,{\frak g}_{-1/2}(U_2)=2(n-2)$,
and ${\frak g}_0(U_2)$ is described by (\ref{iden}). It is clear from
(\ref{iden}) that $\hbox{dim}\,{\frak g}_0(U_2)=n^2-4n+5$. Further,
${\frak g}_{1/2}(U_2)$ and ${\frak g}_1(U_2)$ admit explicit descriptions
(see \cite{S}). These descriptions imply that
$$
{\frak g}_{1/2}(U_2)=\{0\},\qquad {\frak g}_1(U_2)=\{0\},
$$
and therefore
$$
\dim{\frak g}(U_2)=n^2-2n+3<n^2+1.
$$

The proposition is proved. \qed
\medskip\\

We will now finish the proof of Theorem \ref{char2}. It follows from
Proposition \ref{Sig} that, if, for every $p\in M$,  $\alpha_p(I_p)$
is $n^2-2n+1$-dimensional, then $M$ has to be equivalent to either
$B^n$ or $B^{n-1}\times\Delta$ which is impossible since the
dimensions of the automorphism groups of these domains are bigger
than $n^2+1$.

We now assume that $M=M_1\cup M_2$, with $M_1, M_2\ne \emptyset$, and, for $p\in M_1$, we suppose that $\alpha_p(I_p)$ is as the group $G$ in {\bf
(iii)} and, for $p\in M_2$, $\alpha_p(I_p)$ is $n^2-2n+1$-dimensional.
As in the cases $n=2,3$ above, this implies that there exists a
subdomain
$\Omega\subset M$ such that $\hbox{Aut}(\Omega)$ contains a
subgroup $H$ of dimension $n^2+1$ that acts transitively on $\Omega$ and
that contains a subgroup isomorphic to $U(1)\times U(n-1)$. It now
follows from the explicit formulas for the automorphism groups of
$B^n$ and $B^{n-1}\times\Delta$ that such a subgroup $H$ in fact does
not exist. This proves Part {\bf (a)} for $n\ge 4$.

The theorem is proved.\qed
\medskip\\

{\obeylines
Centre for Mathematics and Its Applications 
The Australian National University 
Canberra, ACT 0200
AUSTRALIA 
E-mail address: Alexander.Isaev@anu.edu.au
\hbox{ \ \ }
\hbox{ \ \ }
Department of Mathematics, Campus Box 1146
Washington University in St.\ Louis 
One Brookings Drive
St.\ Louis, Missouri 63130
USA
E-mail address:  sk@math.wustl.edu
}


\begin{thebibliography}{ABCDE}

\bibitem[BDK]{BDK}  Bland, J., Duchamp, T., and Kalka, M., A
characterization of $\CC\PP^n$ by its automorphism group, {\it Complex
Analysis} (University Park, Pa, 1986), 60--65, {\it Lecture Notes In
Mathematics} 1268, Springer-Verlag, 1987.

\bibitem[E] {E} Eisenman, D. A., Holomorphic mappings 
into tight manifolds, {\it Bull.\ Amer.\ Math.\ Soc.} 76(1970),
46--48.

\bibitem[FK]{FK}  Farkas, H.\ and Kra, I., {\it Riemann Surfaces},
Second Edition, Springer-Verlag, 1992.

\bibitem[GIK] {GIK} Gifford, J. A., Isaev, A. V. and Krantz, S. G., On
the dimensions of the automorphism groups of hyperbolic Reinhardt
domains, {\it Illinois J. Math.}, to appear.

\bibitem[GG] {GG} Goto, M.\ and Grosshans, F., {\it Semisimple Lie
 Algebras}, Marcel Dekker, 1978.

\bibitem[GK]{GK} Greene, R.\ E.\ and Krantz, S.\ G.,  
Characterization of complex manifolds by the isotropy subgroups of
their automorphism groups, {\it Indiana Univ.\ Math.\ J.} 34(1985), 865--879.

\bibitem[IK]{IK} Isaev, A. V. and Krantz, S. G., Characterization of Reinhardt domains by their automorphism
groups, Proc. of the 3rd Korean Several Complex Variables Symposium,
Seoul 16-19 December, 1998; to appear in {\it J. Korean Math. Soc.}

\bibitem[Ka]{Ka} Kaup, W., Reele Transformationsgruppen und invariante
Metriken auf komplexen R\"aumen, {\it Invent.\ Math.} 3(1967), 43--70.

\bibitem[KMO]{KMO} Kaup, W., Matsushima, Y. and Ochiai, T., On the
  automorphisms and equivalences of generalized Siegel domains, {\it
  Amer. J. Math.} 92(1970), 475--497.

\bibitem[Ki] {Ki} Kiernan, P., On the 
relations between taut, tight and hyperbolic 
manifolds, {\it Bull.\ Amer.\ Math.\ Soc.} 76(1970), 49--51.

\bibitem[Ko1]{Ko1}   Kobayashi, S.,
{\it Hyperbolic Manifolds and Holomorphic Mappings}, Marcel Dekker,
New York, 1970.

\bibitem[Ko2]{Ko2}  Kobayashi, S., {\it Hyperbolic Complex Spaces},
Springer-Verlag, Berlin, 1998.

\bibitem[N]{N} Nakajima, K., Homogeneous hyperbolic manifolds and
  homogeneous Siegel domains, {\it J.\ Math.\ Kyoto Univ.} 25(1985),
  269--291.

\bibitem[P-S]{P-S} Pyatetskii-Shapiro, I., {\it Automorphic
Functions and the Geometry of Classical Domains} (translated from
Russian), Gordon and Breach,
1969.

\bibitem[S]{S} Satake, I., {\it Algebraic Structures of Symmetric
    Domains}, Kan\^o Memorial Lectures 4, Princeton University Press,
    1980.

\bibitem[U]{U} Upmeier, H. {\it Symmetric Banach Manifolds and Jordan
    $C^*$-Algebras}, North-Holland, 1985.

\bibitem[V]{V} Vigu\'e, J.-P., Les automorphismes analytiques
  isom\'etriques d'une vari\'et\'e complexe norm\'ee, {\it
  Bull. Soc. Math. France} 110(1982), 49--73.

\bibitem[VO]{VO} Vinberg, E.\ and Onishchik, A., {\it Lie Groups and
 Algebraic Groups}, Springer-Verlag, 1990.

\bibitem[Wa]{Wa} Warner, F.\ W., {\it Foundations of 
Differential Manifolds and Lie Groups}, Scott, Foresman \& Co.,
Glenview, London, 1971.

\end{thebibliography}
\end{document}